
\input amstex
\newcount\firstpage\newcount\lastpage\newtoks\runningtitle\newtoks\writer
\firstpage=1
\lastpage=29
\runningtitle={Actions of groups of birationally extendible automorphisms}
\writer={ALAN HUCKLEBERRY AND DMITRI ZAITSEV} 
\def\Bbb{\bold}
\def\C{\Bbb C} 
\def\R{\Bbb R} 
\def\P{\Bbb P} 
\def\Z{\Bbb Z}
\def\N{\Bbb N}

\def\B{\Bbb B}

\def\mapright#1{{\buildrel #1 \over \longrightarrow}}

\def\mapdown#1{\Big\downarrow \rlap{$\scriptstyle#1$}}

\def\o{O}

\def\c{C}
\def\cal{}

\mag=\magstep1
\voffset8pt
\catcode`\@=11
\output={\output@}
\def\output@{\plainoutput}
\hsize=12cm
\vsize=16.1cm
\delimitershortfall=6pt
\nulldelimiterspace=1.44pt
\parindent=18pt
\topskip=12pt
\splittopskip=12pt
\font\titlefont=cmbx10 scaled \magstep1
\font\eightrm=cmr8
\font\eightsl=cmsl8
\font\tensmc=cmcsc10
\def\smc{\relax \ifmmode \Err@ \else \tensmc \fi}

\def\raggedcenter@{\leftskip\z@ plus.4\hsize \rightskip\leftskip
 \parfillskip\z@ \parindent\z@ \spaceskip.3333em \xspaceskip.5em
 \pretolerance9999\tolerance9999 \exhyphenpenalty\@M
 \hyphenpenalty\@M \let\\\linebreak}
\def\uppercasetext@#1{%
   {\spaceskip1.3\fontdimen2\the\font plus1.3\fontdimen3\the\font
    \def\ss{SS}\let\i=I\let\j=J\let\ae\AE\let\oe\OE
    \let\o\O\let\aa\AA\let\l\L
    {\titlefont #1} }}
\def\topmatter{\vskip9pt plus 2pt minus 2pt}
\def\title#1\endtitle{\begingroup\raggedcenter@
   \baselineskip1.3\baselineskip
   \uppercasetext@{\titlefont#1}\endgraf\endgroup}
\def\author#1\endauthor{\vskip27pt\begingroup\raggedcenter@
   \rm \uppercasetext@{\rm #1}\endgraf\endgroup \vskip4pt}
\def\address#1\endaddress{\begingroup\raggedcenter@
   \it#1\endgraf\endgroup}
\def\secondauthor#1\endauthor{
  \vskip9pt
  {\tenpoint\rm\raggedcenter@ and\endgraf}
  \vskip9pt
 \begingroup\raggedcenter@
 \rm\uppercasetext@{#1}\endgraf\endgroup}
\def\abstract#1\endabstract{
  \vskip24pt
  \begingroup
  \raggedcenter@ ABSTRACT\endgraf\endgroup
  \vbox{$$\vbox\bgroup\advance\hsize-6pc\tenpoint\noindent
#1\par\unskip\egroup$$}
  }
\def\endtopmatter{\vskip10pt}
%
%
\newskip\abovesectionskip \abovesectionskip=18pt plus 4pt minus 4pt
\newskip\belowsectiontitleskip \belowsectiontitleskip=12pt plus 2pt minus2pt
\def\penaltyandskip@#1{}
\outer\def\newsection#1\par{
{\penaltyandskip@{-200}\vskip\abovesectionskip
{\bf#1}
\vskip12pt plus 2pt minus2pt}}
\newskip\abovesubsectionskip \abovesubsectionskip=9pt plus 2pt minus 2pt
\outer\def\subsection#1\par{\penaltyandskip@{-100}\abovesubsectionskip
  \leftline{\it#1}\nobreak}
%
\long\def\heading{
{\parindent0pt
\eightrm
\baselineskip6pt
Geometric Complex Analysis \par
edited by Junjiro Noguchi {\eightsl et al.}\par
World Scientific, Singapore, 1995\par
pp.\the\firstpage -- \the\lastpage
\vskip32pt
}}
\headline{\eightrm\ifnum\pageno=\firstpage\hfill\else\ifodd\pageno\hfill\the\runningtitle\hfill\folio\else\folio\hfill{\eightrm\the\writer}\hfill\fi\fi}
\footline{\eightrm\ifnum\pageno=\firstpage\hfill\folio\hfill\else\hfill\fi}

\long\outer\def\claim #1#2\endclaim
    {\medbreak\indent{\bf#1.\enspace}{\sl#2}\par
  \ifdim\lastskip<\medskipamount \removelastskip\penalty55\medskip\fi}
\def\endclaim{}

\long\def\references#1\endreferences{%
\vskip24pt plus 4pt minus 4pt
\centerline{\bf References}
\vskip12pt plus 2pt minus2pt
\leftskip24pt\parindent-24pt#1}
\def\endreferences{}
\long\def\refer[#1]#2\par%
{\hskip0pt\hbox to 16pt{[\hfil#1\hfil]}\hskip8pt#2\par}
\catcode`\@=\active
\pageno=\firstpage

\document
\heading  

\topmatter
   \title
   Actions of groups of birationally extendible automorphisms
   \endtitle
   \author\the\writer    
   \endauthor
  \address  Fakult\"at f\"ur Mathematik \\
	    Ruhr-Universit\"at Bochum \\
	    D-44780 Bochum\\
	    Germany
 \endaddress
\endtopmatter

\newsection 1. Introduction \par
The origin of this work is found in the study of automorphisms of domains
$D$ in $\Bbb C^n$, $n>1$.  For example, suppose for the moment that $D$
is relatively compact and recall that in this case
the group $Aut(D)$ of all holomorphic automorphisms
is a Lie group acting properly on $D$ in the compact-open topology
([6], see also [19]).  It is important to underline the fact that this group is
totally
real so that, compared to holomorphic actions of complex Lie groups,
there is a lack of naturality.

The actions of compact subgroups $K\subset Aut(D)$  extend to the holomorphic
actions of their complexifications.  For example, consider the action of
$K=S^1$ on
an annulus $D$ in the complex plane.  A holomorphic function
$f\in {\cal O}(D)$
has a Fourier (Laurent) series expansion with respect to this action.
This can be regarded as a formal series on $D^\C =\C^*$, where the
complexification $K^\C =\C ^*$ acts holomorphically.
In fact, as a special case of Heinzner's
Complexification Theorem ([12],see also [13]), any domain $D$ equipped with a
compact
group $K$ of holomorphic transformations is naturally contained as a
$K$-stable domain in a Stein manifold $D^\C $ where the reductive group
$K^\C $ acts holomorphically.  If $D$ is Stein, then  it
is just a domain of
convergence for some ``Fourier series'' in $D^\C$.  Thus, except for
convergence questions, in the case of compact groups we are really
confronted with actions of reductive groups.  In this case
the theory of algebraic transformation groups provides us with very stong
tools.

For a non-compact subgroup $G\subset Aut(D)$ the situation is
substantially different.  First of all, as is seen in the simplest example
of the disk $D$ in the complex plane, it is {\it rational} functions which
play an important role.
Secondly, since orbits are non-compact, one is led to study the action
near the boundary.  Without loss of generality we may assume that $G$ is
closed in $Aut(D)$ so that it acts properly and let $p\in \partial D$
be in the closure of some orbit $z\in D$.
The geometry of the action near $p$ is extremely rich.
In fact, under reasonable regularity assumptions,
it might happen that knowledge of the local action near $p$
determines $D$ itself
and general classification results can be proved.
There are numerous indications  of this (see e.g. [8,17,21,22,31,32])
with Rosay's Theorem being the easiest to  state:
{\it If $p$ is a strongly pseudoconvex boundary point,
then $D$ is biholomorphically equivalent to the unit ball
$\B_n:= \{ \sum |z_i|^2 < 1 \}$. } Under far weaker assumptions scaling
methods yield a local description of $D$ near $p$ as also being defined
by polynomial inequalities.

	It is therefore reasonable to begin the study of $G$-actions
on domains by considering the case were $D$ is defined by
polynomial inequalities. In this case, if there is a smooth boundary
point where the Levi form is non-degenerate, by combining results of
Diederich-Pinchuk ([7]) with those in [33],
it follows that $G= Aut(D)$ is a Nash group and its action on $D$
is compatible with the Nash structure. Thus we find ourselves
in the setting of real algebraic geometry.

	Our main results are stated in sections $2-3$.
However, before going to this,
we would like to underline some essential points.
In general, suppose that $G$ acts effectively by
holomorphic transformations on $D$ which extend to rational
transformations of the ambient projective variety $V\subset\P^n$
(e.g.,  if $D$  satisfies Webster's condition $(W)$ below).
The graph of every such transformation defines an $n$-dimensional
cycle in $\P_n\times\P_n$ for some integer $N$.
In this way we obtain a set-theoretic
embedding of $G$ in
the Chow scheme $C_n$ of $n$-dimensional cycles in $\P_n\times\P_n$.
Under certain conditions, which are made precise in the sequel, we show that
$G$ lies in finitely many components of $C_n$. The group
operation on $G$ extends rationally to its Zariski closure $Q$
in $C_n$ and endows $Q$ with a structure of a pre-group in sense of
A.~Weil, which is not a group in general.
The action $G\times D\to D$ extends also
to a rational action $Q\times V \to V$. Again, this is a
pre-transformation space in sense of Weil
which, in general, is not a transformation space.

	Using basic techniques of Weil, we regularize the ``action''
$Q\times V\to V$, i.e. construct an algebraic group $\tilde G$
and an algebraic variety $X$, birationally equivalent to $Q$ and $V$
respectively, such that the induced action $\tilde G\times X\to X$
is regular (Theorems~1.1, 1.3 and 1.4).
As a consequence the action of $G$ extends to a global holomorphic action
of the universal complexification $G^{\C}$ on $X$.

	Further, we employ a ``lifting procedure''
to show that $\tilde G$ is a linear algebraic group. Then a result
of Sumihiro ([16,27])
yields an equivariant embedding of $D$ in a projective
space with a linear action of $\tilde G$ (Theorem~2).
In case of Siegel domains, such equivariant embeddings were obtained by
W.~Kaup, Y.~Matsushima and T.~Ochiai ([14], Theorem~9).

	In the case $D$ is {\it contractible and homogeneous}
under the real analytic action of a connected Lie group $G$ of
birationally extedible automorphisms, R. Penny ([20]) has shown that
the $G$-action extends to a rational action of a real algebraic
group on the ambient space $\C^n$. This is a special case of
Theorem~1 below.

	We would like to conclude this introduction with an
application concerning $G$-invariant meromorphic functions on $D$.
For $x\in D$, let $d(x)$ denote the codimension of
$T_x Gx + iT_x Gx$ in $T_x D$ and $d:= \mathop{max}\limits_{x\in D} d(x)$.
If $f_1,\ldots,f_m$ are $G$-invariant analytically independent meromorphic
functions, then clearly $m\le d$. Now if $\tilde G$ exists as above
and $G$ is Nash, e.g. under the conditions of Corollary~2,  then the
bound $d$ is realized. This follows by applying Rosenlicht's quotient
Theorem ([23]) to the $\tilde G$-action on $X$.

	As indicated  above , we draw our methods from
cycle space theory and algebraic group actions.
On the other hand, our main motivation is of a complex analytic
or representation theoretic nature. Thus we have included details
of results which might be standard in one subject and not so
well-known in the other.

\newsection 2. Algebraic extensions \par

	Here we establish conditions for the existence
of the above extensions which are
birationally equivalent to the ambient space $V\supset D$.
A topological group will always assumed to have a countable basis
at every point.

\medskip
{\bf Definition~2.1}
	An {\bf algebraic extension} of a topological group $G$ of
holomorphic transformations of a domain $D\subset V$ consists
of a homomorphism from $G$ into a complex algebraic group $\tilde G$,
an algebraic variety $X$ birationally equivalent to $V$ via
$\psi\colon V\to X$, such that $\psi|_D$ is a biregular embedding,
and the extension of the action of $G$
to a regular action $\tilde G\times X\to X$.
\medskip

	The existence of algebraic extensions implies, in particular,
that the automorphisms of $D$ which are elements of $G$ extend
to birational mappings from $V$ into itself.
In this case we say that $G$ is a
{\bf group of birationally extendible automorphisms}.

	One condition for the existence of algebraic extensions is given by
the following result.

\claim {Theorem 1}
  Let $V$ be a projective variety,
$D\subset V$ an open set and
$G$ a Lie group of birationally extendible automorphisms of $D$.
Suppose that $G$ has finitely many connected components.
Then there exists an algebraic extension of $G$.
\endclaim

	The existence of an algebraic extension is also equivalent to the
existence of a projective linearization in the following sense.

\medskip
{\bf Definition~2.2}
	A {\bf projective linearization} of a topological group $G$ of
holomorphic transformations of an open set $D\subset V$ consists of
a (continuous) linear representation of $G$ on some $\C^{N+1}$ and
a birational (onto the image) mapping $i\colon V\to\P_N$ such that
the restriction $i|_D$ is biholomorphic and $G$-equivariant.
\medskip

{\bf Remark.} By a rational mapping between two algebraic varieties
$V_1$ and $V_2$ we mean a morphism from a Zariski open dense subset
$U\subset V_1$ into $V_2$. The image is defined to be the (Zariski)
closure of the image of $U$. In general, a point $x\in V_1\setminus U$
may  not correspond to a point of $V_2$.

\claim {Theorem~2}
	Let $V$ be a rational
(i.e. birationally equivalent to $\P_n$) projective variety,
$D$ an open subset of the regular locus of $V$ and
$G$ a topological group of birationally extendible automorphisms of $D$.
Then $G$ has an algebraic extension if and only if it has
a projective linearization.
\endclaim

{\bf Remark.} The condition of rationality of $V$ is perhaps too strong.
However, some condition is needed. If e.g. $D=V=G$
are elliptic curves and $G$ acts on $V$ by translations, this action
coincides with its algebraic extension but has
no projective linearization, because $G$ is a compact complex group.

	Another sufficient condition for the existence of algebraic extensions
is the boundness of the {\bf degree} of
the automorphisms defined by elements of $G$, i.e. the degree of
the graphs $Z_f\subset V\times V$ of the corresponding birational
automorphisms with respect to fixed embedding
$\nu\colon V\times V\hookrightarrow \P_k$.
We also identify $V\times V$ with its image in $\P_k$.

	The boundness of the degree means that the graphs lie as cycles
in finitely many components of the cycle space $\c(V\times V)$.
Therefore, the condition of boundness is independent of
the choice of the embedding $\nu$.

\claim {Theorem~3}
	Let $V$ be a projective variety,
$D\subset V$ an open subset and
$G$ a topological group of birationally extendible automorphisms of $D$.
Then $G$ has an algebraic extension if and only if
the degree of the automorphisms $\phi_g\colon D\to D$
defined by $g\in G$ is bounded.
\endclaim

	In the proof we proceed as follows. Theorem~2 is proven
in section~4. If $\tilde G$ is the algebraic extension, the
rationality of $V$ is used to show that $\tilde G$ is a linear algebraic
group. Then the linearization follows from a theorem of Sumihiro
([16]). The converse in Theorem~2 is straightforward.

	Section~5 is devoted to the proof of Theorem~3.
There we exploit the idea that an action of $\tilde G$ by rational
automorphisms
on $D\subset V$ induces an (almost everywhere defined) mapping $\phi_{\tilde
G}$
from $\tilde G$ into the cycle space $\c(V\times V)$
which can be regarded as a subvariety of the Chow scheme
of an ambient projective space $\P_k$. Here we use the
universal property of the cycle space
([1], see also [5], Proposition~2.20).
The mapping $\phi_{\tilde G}\colon \tilde G\to {C}_n(\P_k)$ is rational
and the boundness of the degree follows from the local constancy of it
on the Chow scheme.

The induced mapping
$\phi_G$ from $G$ into the cycle space $\c(V\times V)$ is continuous
only on an open dense subset $U\subset G$,
but the group operation of $G$
extends to a rational ``group operation'' on the Zariski closure
of $\phi_G(U)$ in $\c(V\times V)$. This operation is defined
via composition of graphs. The objects with rational ``group operations''
were introduced by Weil ([30]) and called {\bf pre-groups}.
The main property is the existence of regularizations of pre-groups,
i.e. algebraic groups which are birationally equivalent to given pre-groups
and the ``group operations'' are compatible with the equivalences.
This property is used to obtain the algebraic group $\tilde G$
for the algebraic extension.

	The next step is to prove that the composition of
$\phi\colon U\to {C}(V\times V)$ and the birational equivalence with $\tilde G$
extends to a continuous homomorphism from $G$ into $\tilde G$.

	The induced ``action'' of $\tilde G$ on $V$ is in general also rational.
Such objects were also introduced
by Weil ([30]) and called {\bf pre-transformation spaces}.
Also in this case he proves the existence of regularizations,
i.e. the (regular) actions of the same group on algebraic varieties
which are birationally equivalent to the original pre-transformation spaces
such that the actions are compatible with the equivalences.
In our case, such regularizations $X$ yield the required algebraic extensions.

	An exposition for pre-groups and pre-transformation spaces
(also not irreducible)
is given in [34].  There, one  also studies the points
where the above regularizations are biregular. This helps
in proving that $D$ is embedded biholomorphically in the context
of Definitions~1.1 and 1.2.

	Theorem~1 is proven in section~6. For this we
use a result of Kazaryan ([15]) to show that the action
$G\times D\to D$ extends to a meromorphic mapping $\tilde G\times V\to V$,
where $\tilde G$ is a complex manifold with $G$ totally really embedded.
Then we prove the boundness of the degree using Proposition~5.1
and the lower semi-continuity of the degree (Lemma~5.1).
Finally, the statement follows from Theorem~3.

\newsection 3. Semialgebraic domains and Nash automorphisms \par

	In general, a domain $D\subset \C^n$ may have no non-trivial
holomorphic automorphisms. On the other hand, in many interesting cases
the automorphism group is very large. The classical examples are
bounded homogeneous domains.
Vinberg, Gindikin and Piatetski-Shapiro ([28]) classified them
and found their canonical realizations as Siegel domains of II kind.
Rothaus ([24]) proved that such realizations are given by
(real) polynomial inequalities. For such reasons, as well as those
mentioned in the introduction, we are intereseted in studiying  domains
defined in this  way. In fact, we
consider more general case of a projective variety $V$ and an open
set $D\subset V$ which is a finite union of the domains given by
finitely many homogeneous polynomial inequalities.
Such set are considered in real algebraic geometry and are called
{\bf semialgebraic} (see e.g. [2] for the elementary introduction to
the theory of semialgebraic sets).

	The following Proposition shows that, for $D$ semialgebraic,
the condition given in Theorem~1 is in some sense also necessary .

\claim {Proposition~3.1}
	Let $V$ be a projective variety,
$D\subset V$ a semialgebraic open subset and
$G$ a topological group of birationally extendible automorphisms of $D$.
Suppose that there exists an algebraic extension of the action of $G$.
Then $G$ is a subgroup of a Lie group $\tilde G$ of birationally extendible
automorphisms of $D$ which extends the action of $G$
to a real analytic action $\tilde G\times D\to D$
and has finitely many connected components.
\endclaim

	Semialgebraic sets are closely related
to the {\bf Nash manifolds} and {\it Nash groups}.
The Nash category is obtained from the real analytic when
we assume all mappings are Nash.
A {\bf Nash mapping} $f\colon D\to D$ is a real analytic mapping,
such that the graph $\Gamma\subset D\times D$ of $f$ is semialgebraic or,
equivalently, Zariski closure $Z_f$ of $\Gamma$
in $V\times V$ has dimension $n=\dim V$.
The reader is refered to [18] and [26] for the precise definitions.

	For a semialgebraic subset $D\subset V$ we prove also the following
criterion.

\claim {Theorem~4}
	Let $V$ be a projective variety,
$D\subset V$ a semialgebraic open subset and
$G$ a topological group of birationally extendible automorphisms of $D$.
The following properties are equivalent:\par
\item {1)}{ $G$ is a subgroup of a Nash group $\tilde G$ of
{\bf birationally extendible automorphisms}
of $D$ which extends the action of $G$ to a Nash action
$\tilde G\times D\to D$;}
\item {2)}{ $G$ is a subgroup of a Nash group $\tilde G$ such that the action
$G\times D\to D$ extends to a Nash action $\tilde G\times D\to D$;}
\item {3)}{ $G$ has an algebraic extension.}
\endclaim

{\bf Remarks.}\par
\item {1)}{ In condition~2 the automorphisms defined by elements
of $\tilde G\setminus G$ are not necessarily birationally extendible.}
\item {2)}{ Since every Nash group is a Lie group with finitely many
components,
Proposition~3.1 is a Corollary of Theorem~4.}

	The proof of Theorem~4 is given in section~7.

	In the remainder of the present
paragraph we mention several applications of Theorem~4
for the bounded semialgebraic domains.
In [33] we gave sufficient conditions on $D$ and $G$
such that $G$ is a Nash group and the action $G\times D\to D$ is Nash.
The domain $D$ is assumed to satisfy the following nondegeneracy condition:

\medskip
{\bf Definition~3.1}
	A boundary of a domain $D\subset\C^n$ is called {\bf Levi nondegenerate}
if it contains a smooth point where the Levi form is nondegenerate.
\medskip

	The group $G$ is taken to be the group $Aut_a(D)$ of all
holomorphic Nash (algebraic) automorphisms of $D$. It was proven in [33]
that, if $D$ is a semialgebraic bounded domain with Levi nondegenerate
boundary,
the group $Aut_a(D)$ is closed in the group $Aut(D)$
of all holomorphic automorphisms and carries
a unique structure of a Nash group such that
the action $Aut_a(D)\times D\to D$ is Nash with respect to this structure.

	Now let $G=Aut_r(D)\subset Aut_a(D)$ be the group of all birationally
extendible
automorphisms of $D$. Then $G$ satisfies the property~2 in Theorem~4
with $\tilde G = Aut_a(D)$. By property~1, $G$ is a subgroup of a Nash
group of birationally extendible automorphisms of $D$.
Since $G$ contains all birationally extendible automorphisms of $D$,
$G$ is itself a Nash group with the Nash action on $D$.
We therefore obtain the following corollary.

\claim {Corollary~1}
	Let $D\subset\C^n$ be a bounded Nash domain with Levi
nondegenerate boundary.
Then the group $Aut_r(D)$ possesses an algebraic extension.
\endclaim

	We now explain sufficient conditions (due to Webster [29])
such that all algebraic automorphisms of $D$ are birationally extendible.
Let $D$ be as in Corollary~1.
The existence of finite stratifications for semialgebraic sets
(see [2], (2.4.4)) implies that the boundary
$\partial D$ is contained in finitely many irreducible real hypersurfaces.
Several of them, let us say $M_1,\ldots,M_k$, have generically
nondegenerate Levi forms. If $\partial D$ is nondegenerate in the sense of
Definition~3.1,
such hypersurfaces exist. The complexification $M^{\C}_i$
of $M_i$ is defined
to be the complex Zariski closures of $M_i$ in $\C^n\times\overline{\C^n}$
where
$M_i$ is embedded as a totally real subvariety
via the diagonal map $z\mapsto (z,\bar z)$.
It follows that $M^{\C}_i$ is an irreducible complex hypersurface.
The Segre varieties $Q_{iw}, w\in \C^n$, associated to $M_i$ are defined by
$$ Q_{iw} := \{z\in \C^n \mid (z, \bar w) \in M^{\C}_i \}.$$
These complexifications and Segre varieties are important biholomorphic
invariants of $D$ and play a decisive role in the reflection
principle which can be used to obtain birational extensions.

\medskip
{\bf Definition~3.2}
	A semialgebraic domain is said to satisfy the condition $(W)$
if, for all $i$, the Segre varieties
$Q_{iw}$ uniquely determine $z\in\C^n$ and $Q_{iw}$ is an irreducible
hypersurface in $\C^n$ for all $z$ in the complement of
a proper subvariety $V_i\subset\C^n$.
\medskip

	A result of Webster ([29], Theorem~3.5) can be formulated
in the following form:

\claim {Theorem~5}
	Let $D\subset\C^n$ be a semialgebraic domain with
Levi nondegenerate boundary
which satisfies the condition $(W)$. Further,
let $f\in Aut(D)$ be an automorphism which is holomorphically
extendible to a smooth boundary point
with nondegenerate Levi form. Then $f$ is birationally extendible to $\C^n$.
\endclaim

{\bf Remark.}
	The mentioned statement of Webster assumes that $f$ extends
biholomorphically to a smooth boundary point where the Levi form is
non-degenerate. By a result of Diederich and Pinchuk ([7]),
this holds for  all automorphisms.

\claim {Corollary~2}
	Let $D\subset\C^n$ be a bounded semialgebraic domain which satisfies
condition $(W)$. Then the whole
group $Aut(D)$ possesses an algebraic extension.
\endclaim

{\bf Acnowledgement.} The authors wish to thank D.~Barlet, K.~Diederich,
 D.~Panyshev and B.~Shiffman for useful discussions.

\newsection 4. Linearization \par

	In the present paragraph we prove Theorem~2.
Assume we are given a projective linearization $i\colon V\to\P_N$.
Let $X$ denote the (Zariski) closure of the (constructible) image $i(V)$.
The subgroup $\tilde G\subset GL_N(\C)$ of all linear automorphisms
of $\P_N$ which preserve $X$ is a complex algebraic subgroup.
Then the pair $(\tilde G,X)$ yields the required algebraic extension.

	The other direction is less trivial. If $G$ has an algebraic
extension, we can assume without loss of generality that $G$ coincides
with the complex algebraic group $\tilde G$.
For the convenience of reader we reformulate here
the conclusion we need to prove.

\claim {Theorem~2'}
	Let $G$ be an complex algebraic group operating regularly on
a rational algebraic variety $X$. Let $D$ be an open set contained
in the regular locus of a quasi-projective subvariety $U\subset X$.
Then there exists a projective linearization.
\endclaim

	Theorem~$2'$ will follow from Lemma~4.1.,
Proposition~4.2. and Sumihiro's Theorem (see below).
Since the regular locus of $X$ is
$G$-invariant, we can replace $X$ with this locus and
Proposition~4.2 can be applied.

\claim {Definition 4.1}
	A line bundle $L$ on an algebraic variety $X$ is called
{\bf birationally very ample} if there exists a finite-dimensional subspace
$W\subset \Gamma(X,L)$ which yields a birational mapping $i_W$ from $X$ into
the corresponding projective space.
\endclaim

\claim {Lemma~4.1}
	Let $G$ be a (complex) algebraic group with a regular action
$\rho\colon G\times X\to X$ on a nonsingular
(not necessarily projective) algebraic variety $X$.
Then there exists a birationally very ample line bundle $L$ on $X$ such that,
for every $g\in G$, $\rho_g^*L\cong L$. If $U\subset X$ is an open
quasi-projective subvariety,
the bundle $L$ and subspace $W\subset \Gamma(X,L)$ can be chosen
such that $i_W$ is regular on $U$.
\endclaim

	{\bf Proof.} Without loss of generality,
$U$ is an open dense quasi-projective subvariety of $X$.
Then the inclusion $\varphi \colon U\to X$ is birational.
Let $C$ be a very ample divisor on $U$ and
$v_0,\ldots,v_N$ a collection of
rational functions on $U$ which yields a basis of ${O}_U(C)$.
The rational functions $v_0\circ\varphi^{-1},\ldots, v_N\circ\varphi^{-1}$
define a birational (onto the image) mapping from $X$ into $\P_N$.
Let $C'$ be the union of polar divisors of all $[\tilde v_i]$ and
$L_{C'}\in Pic(X)$ be the corresponding line bundle. Then $\tilde v_i$'s
can be regarded as sections in $L_{C'}$ which is therefore
birationally very ample.

	It remains to obtain the property $\rho_g^*L\cong L$. The birational
mapping $\varphi\colon \P_n \to X$ is, by definition, a biregular
mapping between Zariski open subsets $U\subset \P_n $ and $U'\subset X$.
Set $E:= \P_n \setminus U$ and $E':=X\setminus U'$ and
let $E_1,\ldots, E_k$ and $E'_1,\ldots, E'_l$ be
the irreducible components of $E$ and $E'$ respectively.
One has the following exact sequences:
$$\oplus_i \Z [E_i] \to Pic( \P_n ) \to Pic(U) \to 0,$$
$$\oplus_i \Z [E'_i] \to Pic(X) \to Pic(U') \to 0.$$
Since $Pic(\P_n)\cong\Z$, it follows that $Pic(U)\cong Pic(U')$ is discrete.
This implies that $Pic(X)$ is discrete. The algebraic group $G$
has finitely many connected components. Therefore, its orbits
in $Pic(X)$ are finite. Thus $G(L_{C'})=\{L_1,\ldots,L_s\}$
as an orbit in $Pic(X)$.
Since the $L_j$'s are birationally very ample, their tensor product
$L:=\otimes_j L_j$ is also birationally very ample and satisfy
the property $\rho_g^*L\cong L$. \hfill {\bf QED}

	We now state and prove a sequence of Lemmas which will yield the
proof of Proposition~4.1.

\claim {Lemma~4.2 }
	Let $G$ and $X$ be arbitrary nonsingular algebraic varieties and $X$ be
birationally equivalent to $\C^n$.
Let $L_{G\times X}$ be a line bundle on $G\times X$.
Then there exist line bundles $L_G$ on $G$ and $L_X$ on $X$ such that
$L_{G\times X}\cong \pi_G^*L_G \otimes \pi_X^*L_X$.
\endclaim

	{\bf Proof.} The special case $X=\C$ is contained in Proposition~6.6.
of Chapter~2 in [11]. By the induction, we obtain the Lemma for $X=\C^n$.
In the general case one has isomorphic Zariski open subsets $U\subset\C^n$
and $U'\subset X$. Set $E:=\C^n\setminus U$ and $E':=X\setminus U'$ and
let $E_1,\ldots, E_k$ and $E'_1,\ldots, E'_l$ be
the irreducible components of $E$ and $E'$ respectively.
Let $L_{G\times X}|G\times U'$ be the restriction and $L_{G\times U}$
its pullback on $G\times U\cong G\times U'$.
Since $L_{G\times U}$ corresponds to a divisor $C$ on $G\times U$,
it is a restriction
of a line bundle $L_{G\times\C^n}$ on $G\times\C^n$ which corresponds to the
closure of $C$ in $G\times\C^n$.
Applying the Lemma to $L_{G\times\C^n}$, we obtain
its splitting  which yields a splitting
$L_{G\times U'}\cong \pi_G^*L_G \otimes \pi_X^*L_{U'}$.
The required splitting for $L_{G\times X}$ is implied now by the surjectivity
of the following map:
$$\oplus_i \Z [G\times E'_i] \oplus Pic(G\times U') \to Pic(G\times X).$$
\hfill {\bf QED}

\claim {Lemma~4.3}
	Let $1\to G_1 \to G \to G_2 \to 1$ be an exact sequence of algebraic
groups. Let either $G$ or both $G_1$ and $G_2$ be linear. Then
all of groups are linear.
\endclaim

	See e.g. [23] for the proof.

\claim {Lemma~4.4}
	Let $G$ be an algebraic group with an effective algebraic action
$\rho\colon G\times X\to X$, where $X$ is a nonsingular
algebraic variety with ${O}^*(X) \cong \C^*$.
Let $L$ be a birationally very ample line bundle on $X$ such that,
for every $g\in G$, $\rho_g^*L\cong L$.
	Then there exists an algebraic
group $\tilde G$ with a surjective homomorphism $\pi\colon \tilde G\to G$
such that \par
\item{1)}{
the action $\tilde G\times X \to X$ defined by $\pi$ is lifted
to an action $\tilde G\times L\to L$,
which preserves the fibres and is linear there;}
\item{2)}{
the kernel of $\pi$ acts effectively  on $L$.}
\endclaim

	{\bf Proof.} Let $\phi\colon \rho^*L \to \pi_G^*L_G \otimes \pi_X^*L_X$
be the isomorphism in Lemma~4.2.
Since $\rho_g^*L\cong L$, one has $L_X\cong L$.

	Let $\tilde G\subset L_G$ be the complement of the zero section.
Our goal now is to define an algebraic group structure on $\tilde G$ and
to construct an algebraic action $\tilde G\times L\to L$.
The action $\tilde\rho\colon \tilde G\times L\to L$
is defined as the composition
$$
L_G\times L \to \pi_G^*L_G \otimes \pi_X^*L \to \rho^*L \to L,
\tag 4.1$$
where the first mapping is given by two isomorphisms
$L_G\times X\to \pi_G^*L_G$ and $G\times L\to \pi_X^*L$.
The composition (4.1)
makes the following diagram commutative:
$$\matrix
\tilde G\times L         & \to & L           \\
\downarrow               &     & \downarrow  \\
G \times X               & \to & X            \endmatrix $$
Let $g\in G$ be fixed. Then the fibre $(L^*_G)_g (\cong \C^*)$ of $\tilde G$
over $g$ defines a $1$-dimensional family of automorphisms of $L$ which
lift the automorphism $\rho_g\colon X\to X$ defined by the action of $G$.
(By an automorphism of $L$ we mean an algebraic isomorphism of $L$ onto
itself which takes fibres in fibres and is linear on them.)
Since ${O}^*(X) \cong \C^*$, such automorphisms of $L$ form a
$1$-dimensional family which coincides therefore with the family defined by
$(L^*_G)_g$. We obtain a one-to-one correspondence between
the elements of $\tilde G$ and the liftings of the automorphisms
$\rho_g\colon X\to X$ for $g\in G$.

	The set of all automorphisms of $L$ which lift $\rho_g$ for some $g\in G$
forms a group in a natural way. The above one-to-one correspondence
transfers this group structure to $\tilde G$. The regular mapping in
(4.1) defines a group action $\tilde G\times L\to L$
with respect to this structure.

	Now we wish to prove that the group operation
$\tilde G\times \tilde G\to \tilde G$, $(g,h)\to gh$ is algebraic.
Since the action $\tilde\rho \colon \tilde G\times L\to L$ is algebraic,
the map
$$
	\alpha\colon \tilde G\times \tilde G\times L \to L,
	(g,h,l)\mapsto \tilde\rho(g,\tilde\rho(h,l))
$$
is also algebraic. We find the product $t:=gh \in \tilde G$ from the relation
$$
\alpha(g,h,l) = \rho (t,l).
\tag 4.2$$

	For a fixed arbitrary point $l_0\in \tilde G$, the mapping
$$
\imath:=\pi\times\rho(\cdot,l_0) \colon \tilde G \to G\times L
$$

is a regular embedding. By (4.2), $t = gh$ can be expressed as follows:
$$
t(g,h) = \imath^{-1} \circ (\pi(g)\pi(h), \alpha(g,h,l_0).
$$

This proves the algebraicity of the group operation on $\tilde G$.
It remains to prove that the inverse map
$\tilde G\to \tilde G, g\mapsto g^{-1}$ is also regular.
For this consider
$$
\Gamma:= \{(g,h)\in \tilde G\times \tilde G \mid
	\tilde\rho(g,\tilde\rho(h,l)) = l \hbox{ for all } l\in L\}.
$$
This is the graph of $g\mapsto g^{-1}$ which projects bijectively
on both factors $\tilde G$. Since $\Gamma$ is an algebraic subset
and $\tilde G$ is nonsingular, the inverse mapping is regular.
\hfill {\bf QED}

	The following is a foundational result for algebraic group actions.

\claim {Lemma~4.5}
	Let $G\times X \to X$ be an algebraic action of an algebraic group
$G$ on an algebraic variety $X$ which lifts to an action on a line bundle
$L$ on $X$. Then the induced action on the space of sections
$\Gamma(X,L)$ is rational and locally finite.
\endclaim

	The proof coincides with the proof of Lemma~2.5. in [16],
where $G$ is regarded as an arbitrary algebraic group.

\claim {Proposition~4.1}
	Let $G$ be an algebraic group with an effective algebraic action
$\rho\colon G\times X\to X$, where $X$ is a rational nonsingular
algebraic variety with
${O}^*(X) \cong \C^*$.
Let $L$ be a birationally very ample line bundle on $X$ such that,
for every $g\in G$, $\rho_g^*L\cong L$. Then $G$ is linear algebraic.
\endclaim

	{\bf Proof.} Let $W\subset\Gamma(X,L)$ be
the finite dimensional subspace in Definition~4.1. By Lemma~4.5.
applied to the group $\tilde G$, $W$ generates a finite dimensional invariant
subspace $\tilde W\subset\Gamma(X,L)$ which also yields a birational
mapping $i_{\tilde W}\colon X\to \P(\tilde W^*)$.
We obtain a representation of $\tilde G$ in $\tilde W^*$.
Let $K\subset \tilde G$ be its kernel. An element $k\in K$ acts
trivially on $i_{\tilde W}(X)$ and therefore on $X$. Since the action
$G\times X\to X$ is effective, this implies that $K\subset Ker\,\pi$.
But the kernel of $\pi$ acts effectively on $\Gamma(X,L)$ which implies
$K=\{e\}$.

	Thus, $\tilde G$ is a linear algebraic group. Since $G$ is a homomorphic
image of $\tilde G$, it is also linear algebraic.
\hfill {\bf QED}

\claim {Lemma~4.6}
	Let $X$ be a rational nonsingular algebraic variety and
$\rho\colon G\times X\to X$ be an algebraic action of an algebraic group $G$
which satisfies the property ${O}^*(G) \cong \C^*$.
Let $L$ be a birationally very ample line bundle on $X$ such that
for every $g\in G$, $\rho_g^*L\cong L$. Then the action of $G$ on $X$
is trivial.
\endclaim

	{\bf Proof.} We prove the Lemma by induction on $\dim X$.
The condition ${O}^*(G) \cong \C^*$ implies the connectedness and
irreducibility of $G$. Let $\dim X = 0$. Then $X$ is discrete
and the action is trivial.

Now assume $\dim X \ge 1$. Let $x_0\in X$ be an arbitrary point and define

$$
	F(x_0):=\{x\in X \mid \forall f\in {O}^*(X), f(x)=f(x_0) \} \subset X.
$$

	Since ${O}^*(G) \cong \C^*$, the orbit $Gx_0$ lies in $F(x_0)$.
This is true for any orbit $Gx$ with $x\in F(x_0)$ and therefore
$F(x_0)$ is $G$-invariant. Let $X_0\subset X$ be an irreducible component
with $x_0\in X_0$. Since $G$ is irreducible, $X_0$ and $F'(x_0):=X_0\cap
F(x_0)$
are also $G$-invariant.

	Now two cases are possible. If $\dim F'(x_0)< \dim X$, the action on $F'(x_0)$
is trivial by induction. If $\dim F'(x_0)= \dim X$, then $F'(x_0)=X_0$
and ${O}^*(X_0)\cong \C^*$. By Proposition~4.1,
$\hat G:=G/Ker(\rho|_{X_0})$ is a linear algebraic group. The condition
${O}^*(G) \cong \C^*$ for $G$ implies the same condition for $\hat G$.
Since $\hat G$ is linear algebraic, it is trivial. Thus,
$Ker(\rho_{X_0})=G$ which means that the action on $X_0$ is trivial.

	In summary we obtain that, for every $x_0\in X$ and $g\in G$, $gx_0=x_0$.
This means that $G$ acts trivially. \hfill {\bf QED}
\medskip

	The following is straightforward.

\claim {Lemma~4.7}
	Let $G$ be an algebraic group and $e\in G$ the unit. Then the subvariety
$$
F(e):=\{g\in G\mid \forall f\in {O}^*(G), f(g)=f(e) \}
$$
is an algebraic subgroup.
\endclaim

\claim {Lemma~4.8}
	Let $G$ be an algebraic group such that the global invertible regular
functions separate its points. Then $G$ is linear algebraic.
\endclaim

	This is a corollary of the following lemma:

\claim {Lemma~4.9}
	Let $G$ be an algebraic group such that the global regular
functions separate points of it. Then $G$ is linear algebraic.
\endclaim

	{\bf Proof.} By Corollary~3. in [23], page 431, there exists
an algebraic subgroup $D\subset G$ such that the quotient $G/D$ is linear
and such that the kernel of any algebraic homomorphism from $G$ into a
linear group contains $D$. It is enough to prove that $D=\{e\}$.

	Assume the contrary. Let $g\ne e$ be an arbitrary point in $D$. Since
the points of $G$ are separated by global regular functions, there exists
a function $f\in {O}(G)$ such that $f(g)\ne f(e)$. By Lemma~4.5,
$f$ generates a finite dimensional $G$-invariant subspace
$W\subset {O}(G)$. The canonical representation of $G$ in $W$
is a homomorphism from $G$ into a linear group such that its kernel
does not contain $g$. This contradicts to the property of $D$ and the
fact that $g\in D$. \hfill {\bf QED}

	Now we drop the assumption ${O}^*(X) \cong \C^*$ in
Proposition~4.1.

\claim {Proposition~4.2}
	Let $G$ be an algebraic group with an effective algebraic action
$\rho\colon G\times X\to X$ on an algebraic variety $X$.
Let $L$ be a birationally very ample line bundle on $X$ such that,
for every $g\in G$, $\rho_g^*L\cong L$. Then $G$ is linear algebraic.
\endclaim

	{\bf Proof.} We proceed by induction on $\dim G$. The Proposition is
trivial for $\dim G=0$.

By Lemma~4.3, we can assume $G$ to be connected.
Let $F(e)$ be the algebraic subgroup defined in Lemma~4.7.
If $\dim F(e)=\dim G$, it follows that $F(e)=G$ which implies
${O}^*(G)\cong \C^*$. Then, by Lemma~4.6, $G$ acts trivially.
Since it acts also effectively it is trivial (and of course linear algebraic).

If $\dim F(e) < \dim G$, the subgroup $F(e)$ is linear by the induction.
The group $G/F(e)$ satisfies conditions of Lemma~4.8 and is
also linear algebraic. Then, by Lemma~4.3, $G$ itself is linear.
\hfill {\bf QED}

\claim {Sumihiro' Theorem}
	Let $G$ be a {\bf linear} algebraic group operating regularly on
a rational algebraic variety $X$. Let $D$ be an open set contained
in the regular locus of a quasi-projective subvariety $U\subset X$.
Then there exists a projective linearization.
\endclaim

	The original proof ([16,27]) is given for the case
$D$ is an orbit or $U$ is $G$-invariant.
In the general case we take a birationally very ample bundle $L$
given by Lemma~4.1 and follow the proof in [16].

\newsection 5. Algebraic extensions for bounded degree \par

	The goal of this section is to prove Theorem~3.

	Recall that by the {\it degree} of $\phi_g$
with respect to a fixed (biregular)
embedding $\nu\colon V \times V \to \P_k$ we mean
the degree of the closed graph $Z_g\subset V \times V$ of $\rho_g$
embedded in $\P_k$ via $\nu$.

	We use the following universal property of the cycle space
([1], see also [5], Proposition~2.20):

\claim {Proposition~5.1 }
	Let $X$ and $S$ be irreducible complex spaces. There exist a natural
identification between:\par
\item {1)}{ meromorphic maps $\phi\colon S\to {C}_n(X)$, and}
\item {2)}{ $S$-proper pure $(d+n)$-dimensional cycles $F$
of $S\times X$ ($d=\dim S$).}
\endclaim

	Let $\tilde G$ be an algebraic extension of $G$ as in
Definition~2.1 and $\Gamma\subset \tilde G\times V\times V$
be the graph of the rational ``action'' of the algebraic group $\tilde G$.
By Proposition~5.1, this action induces a rational
mapping $\mu\colon \tilde G\to {C}_n(\P_k)$. The finiteness of the
number of irreducible components of $\tilde G$ implies the boundness
of the degree of $\phi_g$ for all $g$ in an open dense subset of $\tilde G$.
The global boundness is obtained by the following lemma.

\claim {Lemma~5.1}
	The degree is a lower-semicontinuous function on $G$.
\endclaim

	{\bf Proof.} Let $g_0\in G$ be an arbitrary point and $g_m,m\in\N$
an arbitrary sequence with $g_m\to g_0$. It is enough to prove that
$deg(Z_{g_m})\ge deg(Z_{g_0})$ up to finite set of $m\in\N$.
Assume on the contrary that $deg(Z_{g_m}) < deg(Z_{g_0})$
for a subsequence which is again denoted by $g_m$.
By a theorem of Bishop ([3]),
$Z_{g_m}$ can be assumed to converge to some
cycle $Z_0$ with $deg(Z_0)<deg(Z_{g_0})$. By the continuity of
the action $\rho$, one has $Z_{g_0}\cap D\times D \subset Z$, which implies
$Z_{g_0} \subset Z_0$. On the other hand,
by the continuity of degree (which is equivalent
to the continuity of the volume), $deg(Z_0)<deg(Z_{g_0})$,
which contradicts the above inclusion. \hfill {\bf QED}

\bigskip
{\it 1. Formulation.}
\medskip

	The other less trivial direction in Theorem~3
will be a corollary of the following statement:

\claim {Theorem~$3'$}
	Let $D\subset V$ be an open subset in a projective variety $V$,
$G$ a topological group and $\rho\colon G\times D\to D$
a continuous action such
that, for every $g\in G$, the homeomorphism $\rho_g\colon D\to D$ extends
to a birational mapping from $V$ into itself (which we also denote by
$\rho_g$).
Assume that the set of degrees of all $\rho_g, g\in G$ is bounded.
Then there exist:\par
\item {1)}{ an algebraic group $\tilde G$,}
\item {2)}{ a continuous homomorphism $\phi\colon G\to \tilde G$,}
\item {3)}{ an algebraic variety $X$,}
\item {4)}{ an algebraic action $\tilde G\times X\to X$,}
\item {5)}{ a birational mapping $\psi\colon V \to X$
	such that $\psi|_D$ is biholomorphic and $G$-equivariant.}
\endclaim

\bigskip
{\it 2. Properties of the group $G$.}
\medskip
We begin by noting an elementary basic fact.
\claim {Lemma~5.2}
	Let $G$ be a topological space and
$f\colon G\to \Z$ a lower-semicontinuous
function which is bounded from above.
Then the set $U\subset G$ of all local maximums of $f$ is
 open and dense in $G$.
Moreover, $f|_U$ is locally constant.
\endclaim

	Let $Z\subset G\times V \times V $ and
$\nu(Z)\subset G\times\P_k$ be the families of all $Z_g$ and $\nu(Z_g)$,
$g\in G$, respectively. We denote by $U\subset G$
the set of all local maxima of the degree,
which is open dense by Lemma~5.2.

\claim {Lemma~5.3}
	Let $U$ be a topological space and $\{\phi_g\}_{g\in U}$
a continuous family of automorphisms of $D$ which exted to birational
mappings from $V$ to $V$ with (closed) graphs $Z_g$.
 Assume that the degree of $Z_g$ is locally constant on $U$.
Let the automorphisms depend continuously on $u\in U$.
Then the family $Z$ is closed in $U\times V\times V$.
\endclaim

	{\bf Proof.} Let $(g_0,z_0)\in U\times V\times V$ be a point and
$(g_m,z_m)\to (g_0,z_0)$ a sequence with $z_m\in Z_{g_m}$.
By a theorem of Bishop ([3]), the sequence of cycles $Z_{g_m}$
can be assumed to converge to some $Z_0$. By the continuity
of the automorphisms,
one has $Z_{g_0} \subset Z_0$.
Since $g_0\in G$ is a local maximum of the degree
and the degree is a continuous function on cycles, one obtains
$deg Z_0 \le deg Z_{g_0}$. This means $Z_{g_0} = Z_0$ and $(g_0,z_0)\in Z$.
\hfill {\bf QED}

	Using the family
$Z\cap (U\times V\times V)\subset U\times\P_k$ we define
a continuous mapping $\phi$ from $U$ into the Chow
scheme $C$ of cycles in $\P_k$ (see [10,25]).
We recall briefly the construction of the components of $C$.
Let $Z_g,g\in P$ be an arbitrary family  of
irreducible subvarieties of $\P_k$ of fixed dimension $n$ and degree $d$,
parameterized by a set $U$. The $n+1$-tuples $(H_0,\ldots,H_k)$ of
hyperplanes in $\P_k$ are parameterized by $S:=(\P_k^*)^{n+1}$.
We define $V_g\subset \P_k\times S$ by
$$V_g:=\{(z,H_0,\ldots,H_n \mid z\in Z_g\cap H_0\cap\cdots\cap H_n,$$
and denote by $\pi(V_g)\subset S$ its projection.
Then all $V_g$'s, $g\in U$ and, therefore, all $\pi(V_g)$'s are irreducible
subvarieties. Moreover, $\pi(V_g)$'s are of codimension $1$
and of multidegree $(d,\ldots,d)$. They are given
uniquely up to multiplications by constants
by multihomogeneous polynomials $R_g\subset \C[S]_{d,\ldots,d}$
of multidegree $(d,\ldots,d)$.

	Let $\P_N= \P(\C[S]_{d,\ldots,d})$ denote the projectivization
of the space of such polynomials and $N_R\subset S,R\in \P_N$
the family of zero sets of them.
Therefore we obtain a mapping $\phi \colon U\to \P_N$ which
associates to every $g\in U$ the Chow coordinates $[R]=[R_g]\in \P_N$
of $Z_g$ such that $\pi(V_g)=N_R$.

We utilize the following topological universal property of the Chow scheme:

\claim {Proposition 5.2}
	Let $U$ be a topological space
and $Z_g\in\P_k,g\in U$, a closed family, i.e. the subset
$$Z=\{ (g,z)\mid z\in Z_g \} \subset U\times \P_k$$ is closed.
Suppose that the dimension and degree of $Z_g$ are constant.
Then $\phi\colon U\to C$ is a continuous mapping.
\endclaim

{\bf Proof.}
The closedness of $Z_g,g\in U$, implies the closedness of
$V_g,g\in U$, because the latter is defined by a closed condition.
Since the projective space $\P_k$ is compact,
the family of projections $\pi(V_g)$ is also closed.
The graph $\Gamma\subset U\times \P_N$ of the mapping $\phi$
is defined by the condition
$$ \Gamma=\{ (g,[R]) \mid Z_g\subset N_R \}. $$
It is sufficient to prove that $\Gamma\subset U\times \P_N$
is closed.

	Let $(g_0,[R_0])$ be a point in the complement of $\Gamma$.
This means that $R_0(z_0)\ne 0$ for some $z_0\in Z_{g_0}$.
Then there exist a neighborhood $U(z_0)\subset\P_k$ of $z_0$
and a neighborhood $U(R_0)\subset \P_N$ of $[R_0]$
such that $R(z)\ne 0$ for all $z\in U(z_0)$ and $[R]\in U(R_0)$.
We claim that there exists a neighborhood $U(g_0)\subset U$
such that $U(z_0)\cap Z_g\ne\emptyset$ for all $g\in U(z_0)$.
Indeed, otherwise there would be a sequence $g_m\to g_0$
without this property for $Z_{g_m}$.
By a theorem of Bishop ([3]),
one has, passing if necessary to a subsequence,
$Z_{g_m}\to Z_{g_0}$, which is a contradiction.

	Therefore, the whole neighborhood $U(p_0)\times U(R_0)$ of $(g_0,[R_0])$
belongs to the complement of $\Gamma$. This proves the closedness
of the graph $\Gamma$ which means the continuity of the mapping $\phi$.
\hfill {\bf QED}

	The Chow scheme $C$ is a collection of projective varieties parameterized
by the dimension and degree of cycles. In Theorem~5.4 we assume
that the set of all degrees of $Z_g,g\in G$ is bounded.
Therefore, the image $\phi(U)$ is contained in finitely many components
of the Chow scheme. Let $Q$ denote the Zariski closure of $\phi(U)$
in $C$. It is a projective variety. Let $F\subset Q\times\P_k$ denote
the universal family over $Q$. Since $F_v\subset  V \times V $ for
all $v$ from the Zariski dense subset $\phi(U)$, one has
$F\subset Q\times V\times V$.

\claim {Lemma~5.4}
	Let $Q$ be an algebraic variety and $F\subset Q\times V\times V$
a closed algebraic family of subvarieties $F_v\subset V\times V,v\in Q$,
of pure dimension $n$.
For every $v$ from a Zariski dense subset $\phi(U)\subset Q$, assume that
the fibre $F_v$ is
the closed graph of a birational mapping $\rho_v\colon V \to V$.
Then this is true for all $v$ from a Zariski {\bf open} dense subset $Q'$
with $\phi(U)\subset Q'\subset Q$. Moreover, there exists a Zariski open
subset $F'\subset F$ which intersects every graph $F_v,v\in Q'$, along
a Zariski dense graph of a {\bf biregular} mapping $\phi'_v$.
\endclaim

	{\bf Proof.} Let $Q_1\subset F$ be the set of all
$(v,x)\in Q\times V$ such that the fibres $F_{(v,x)}\subset V$
are finite. Since the fibre dimension is upper-semicontinuous, $Q_1$
is a Zariski open subset of $Q\times V$. The family $F$ is a finite
ramified covering of $Q_1$.
The set
$$R:= \{(v,x)\in Q\times V \mid \rho_v \hbox{ is biregular at $x$ }\}.$$
is a dense subset of $Q_1$ and the fibres $F_{(v,x)}$ over $R$ consist
of single points.  Therefore the covering $F$ has only one sheet and every
fibre $F_{(v,x)}$, $(v,x)\in Q_1$, consists of a single point.
If, for some $v\in Q$, $(\{v\}\times V)\cap Q_1$ is dense in $(\{v\}\times V)$,
this means that $F_v\subset V\times V$ is the graph of a rational mapping
$\rho_v\colon V\to V$. This is true for all $v$ from a Zariski open dense
subset
$Q'_1, \phi(U)\subset Q'_1\subset Q$, which can be taken to be the intersection
of the projections of irreducible components of $Q_1$ on $Q$.

	Similarly, using the projection on the product of $Q$ and the other
copy of $V$, we can construct Zariski open dense subsets
$Q_2\subset Q\times V$, $Q'_2, \phi(U)\subset Q'_2\subset Q$,
such that $\phi_v^{-1}$ is regular at $x\in V$ for all $(v,x)\in Q_2$
and $F_v$ is a graph of a mapping $\phi_v$ with rational inverse for
all $v\in Q'_2$. Then the intersection $Q':=Q'_1\cap Q'_2$
satisfies the required properties.

	The required Zariski open subset $F'\subset F$ can be given by the
formula
$$F':=\pi_1^{-1}(Q_1)\cap \pi_2^{-1}(Q_2),$$
where $\pi_1,\pi_2\colon F\to Q\times V$ denote the projection
on the product of $Q$ and the first (resp. the second) copy of $V$.
\hfill {\bf QED}

	Now we wish to extend the group operation
$U\times U\to G$ to a rational mapping $Q\times Q\to Q$.
Let $Q'\subset Q$ be given by Lemma~5.4.

\claim {Lemma~5.5}
	There exist rational mappings
$\alpha,\alpha_1,\alpha_2\colon Q\times Q\to C$
such that for all $(v,w)\in Q'\times Q'$, where $\alpha$
(resp. $\alpha_1$ and $\alpha_2$)
is defined, the fibre $F_{\alpha(v,w)}$
(resp. $F_{\alpha_1(v,w)}$ and $F_{\alpha_1(v,w)}$)
coincides with the closed graph of the birational correspondence
$\rho_v\circ\rho_w$
(resp. $\rho_v\circ\rho_w^{-1}$ and $\rho_w^{-1}\circ\rho_v$).
\endclaim

	In the construction of the mappings $\alpha$, $\alpha_1$ and $\alpha_2$
we use the following algebraic universal property of the Chow scheme.
Recall that $C$ denotes the Chow scheme of $\P_k$ and
$F\subset C\times\P_k$ the universal family over $C$.

\claim {Proposition~5.3}
	Let $X$ be a quasi-projective variety,
$Z\subset X\times\P_k$ a closed pure-codimensional subvariety.
Then there exists a rational mapping
$i\colon X\to C$ with $Z_g=F_{i(v)}$ for all $v\in X$,
such that $i$ is regular at $v$.
\endclaim

	This is a con\-sequence of Propo\-sition~5.1 and Chow's
theorem ([9], p.~167).

{\bf Proof of Lemma~5.5.}
	We construct here the extension $\alpha_1$ of the mapping
$(g,h)\mapsto gh^{-1}$. The construction of $\alpha$ and $\alpha_2$
is completely analogous.

	The idea of construction is to consider the family of graphs
of $gh^{-1}\colon V \to V $ over $Q\times Q$ and to utilize
the above universal property for it.
Let $W_1$, $W_2$ and $W_3$ denote different copies of $V$
and $\pi_1$ ,$\pi_2$ and $\pi_3$ be the projections of
$W_1\times W_2\times W_3$ onto $W_2\times W_3$, $W_1\times W_3$
and $W_1\times W_2$ respectively. Then, for $g_1,g_2\in U$,
the graph of $\phi_{g_1g_2^{-1}}$ is equal to the closure of
$$
Z'_{g_1g_2^{-1}} = \pi_2(\pi_3^{-1}(Z'_{g_2})\cap\pi_1^{-1}(Z'_{g_1})),
\tag 5.1$$
where  $Z'_{g_1}\subset W_2\times W_3$ and
$Z'_{g_2}\subset W_2\times W_1$ are the regular parts of the graphs of
$\rho_{g_1}\colon W_2\to W_3$ and
$\rho_{g_2}\colon W_2\to W_1$ respectively.

Using formula (5.1) we define a constructible family
$\tilde F\subset Q'\times Q'\times W_1\times W_3$:
$$
\tilde F = \pi_2(\pi_3^{-1}F'_2\cap\pi_1^{-1}(F'_1)),
\tag 5.2$$
where $F'_1\subset Q\times W_2\times W_3$ and
$F'_2\subset Q\times W_2\times W_1$ are different copies
of $F'\subset Q\times V \times V $. This is given by Lemma~5.4

	By the choice of $F'$ and $Q'$, every fibre $\tilde F_{v_1,v_2}$,
$v_1,v_2\in Q'$ is purely $n$-dimensional. Therefore the
family $\tilde F$ is closed and of locally constant degree
in a Zariski open dense subset
$Q''\subset Q'\times Q'$. By Proposition~5.3,
there exists a rational mapping $\alpha_1\colon Q''\to C$ with
$\tilde F_{(v_1,v_2)}=F_{\alpha_1(g)}$.
Since $Q''$ is Zariski open and dense,
$\alpha_1$ extends to a rational mapping $Q\times Q\to C$
which has the required properties. \hfill {\bf QED}

	Since the maps $\alpha$ and $\alpha_1$ extend the group operations,
we write $\alpha(v,w) = vw$ and $\alpha_1(v,w) = vw^{-1}$ whenever these
values are defined.

\claim {Lemma~5.6}
	The mapping $(v,w)\mapsto (vw,w)$ is injective on $Q'\times Q'$.
\endclaim

	{\bf Proof.} Let $v,w\in Q'$ be arbitrary points. By Lemma~5.4,
the fibres $F_v,F_w$ are closed graphs of birational mappings
$\rho_v,\rho_w\colon V \times V $. By Lemma~5.5,
the fibre $F_{vw}$ is the
closed graph of the composition $\rho_v\circ\rho_w$.
If $v_1w=v_2w$, their fibres are also equal which implies the equality
$\rho_{v_1}\circ\rho_w=\rho_{v_2}\circ\rho_w$.
Since $\rho_w\colon V \to V $ is birational, we obtain
$\rho_{v_1}=\rho_{v_2}$, which means $v_1=v_2$. \hfill {\bf QED}

	The following Lemma states the existence of right and left divisions
of ``generic'' elements.

\claim {Lemma~5.7}
	The mappings $(v,w)\mapsto (vw,w)$ and $(v,w)\mapsto (wv,w)$
are birational mappings
from $Q\times Q$ into itself with the inverses
$(v,w)\mapsto (vw^{-1},w)$ and $(v,w)\mapsto (w^{-1}v,w)$.
The variety $Q$ is pure-dimensional.
\endclaim

	{\bf Proof.} We prove the statement for the first mapping. The proof
for the second one is completely analogous.
We first wish to prove that the closed image
of $Q\times Q$ under the mapping $(v,w)\mapsto (vw,w)$ lies in $Q\times Q$.
Since $\phi(U)$ is Zariski dense in $Q$,
$\phi(U)\times\phi(U)$ is Zariski dense in $Q\times Q$.
Then the subset $W$ of $\phi(U)\times\phi(U)$, where $vw$ is defined,
is also Zariski dense. Since the mapping $G\times G\to G\times G$,
$(g,h)\to (gh,h)$ is a homeomorphism, the preimage $U'$ of $U\times U$
is open dense in $G\times G$ and therefore $U'':=U'\cap \phi^{-1}(W)$ is open
dense in $\phi^{-1}(W)$. This implies that $Q'':=(\phi\times\phi)(U'')$
is Zariski dense in $W$ and thus in $Q\times Q$.

	Let $(v,w)=(\phi(g),\phi(h))\subset Q''$ be an arbitrary point.
By Lemma~5.5,
the fibre $F_{vw}$ is the closed graph of the composition $\rho_v\circ\rho_w$.
The latter birational mapping coincides with the automorphism
$\rho_{gh}$ defined by $gh\in U$. This means that $vw\in Q$.
Since $Q''$ is Zariski dense, this inclusion is valid for all
$(v,w)\in Q\times Q$ where $vw$ is defined. Thus the mapping
$(v,w)\mapsto (vw,w)$ is a rational mapping from $Q\times Q$ into itself.

	The projective variety $Q$ has finitely
many irreducible components. Let $Q_0\subset Q$ be a component
of maximal dimension and $Q_1$ an arbitrary component.
Then $(vw,w)\in Q_2\times Q_1$ for $(v,w)\in Q_0\times Q_1$,
where $Q_2$ is also a component of $Q$. By the choice of $Q_0$,
$\dim Q_2\le\dim Q_0$. Since the restriction on $Q_0\subset Q_1$
of the mapping in Lemma~5.6 is injective on the open subset
$(Q_0\times Q_1)\cap(Q'\times Q')$, its closed image coincides
with $Q_2\times Q_1$. Therefore the composition of
$(v,w)\mapsto (vw,w)$ and $(v,w)\mapsto (vw^{-1},w)$ is defined in
an open dense subset $Q''$ of $Q_0\times Q_1$. It is equal to the identity on
the Zariski dense subset $Q''\times \phi(U)$, i.e. it is the identity.
By the injectivity in Lemma~5.6,
the mapping $(v,w)\mapsto (vw,w)$ is birational
from $Q_0\times Q_1$ into $Q_2\times Q_1$ with the
inverse $(v,w)\mapsto (vw^{-1},w)$.
In particular, $\dim Q_2 =\dim Q_1$.

	Now, by Lemma~5.6, the components $Q_2$ are different for
different $Q_1$ and fixed $Q_0$. If $Q_1$ runs through all components,
$Q_2$ also does. This implies that $Q$ is pure-dimensional.
Thus, we can take for $Q_0$ and $Q_1$ any two components and repeat
the above proof. \hfill {\bf QED}

\claim {Lemma~5.8}
	Let $Q'\subset Q$ be as in Lemma~5.4 and $vw$ be defined
and in $Q'$ for $v,w\in Q'$. Then $\rho_{vw}=\rho_v\circ\rho_w$.
\endclaim

	{\bf Proof.} In case $v=\phi(g),w=\phi(h)$ for $g,h,gh\in U$
one has $\rho_{vw}=\rho_{gh}=\rho_g\circ\rho_h=\rho_v\circ\rho_w$.
Since the set of above points $(v,w)\in Q'\times Q'$ is Zariski dense,
the required relation is valid in general.\hfill {\bf QED}

	We now establish the associativity of the operation
$(v,w)\to vw$.

\claim {Lemma~5.9}
	Let $u,v,w\in Q$ be arbitrary points. Then $(uv)w=u(vw)$ whenever
both expressions are defined.
\endclaim

	{\bf Proof.} By Lemma~5.7,
the above expressions are defined on a Zariski open dense subset
$Q''\subset Q^3$. For $u,v,uv,vw \in Q'$ the fibres of both expressions
are the graphs of $\rho_u\circ\rho_v\circ\rho_w$ by Lemma~5.8.
The latter set is Zariski dense. \hfill {\bf QED}

\claim {Lemma~5.10}
	The operation $(v,w)\to vw$ induces a group structure on the set $S$
of all irreducible components of $Q$.
\endclaim

	{\bf Proof.} The associative property follows from Lemma~5.9.
Lemma~5.7 implies the existence of right division in $S$. The
existence of left division is proved analogously by using the
birational correspondence $(v,w)\to (v,v^{-1}w)$. \hfill {\bf QED}

	It follows from Lemmas~5.7 and 5.9 that $Q$ is an
{\bf algebraic pre-group} in sense of [34].
Recall that an {\bf algebraic pre-group}
is an algebraic variety $V$ with a rational
mapping $V\times V\to V$, written as $(v,w)\mapsto vw$, such that:\par
\item {1)}{
for generic $(u,v,w)\in V\times V\times V$ both expressions
$(uv)w$ and $u(vw)$ are defined and equal
(generic associativity condition);}
\item {2)}{
the mappings $(v,w)\mapsto (v,vw)$ and $(v,w)\mapsto (v,wv)$
from $V\times V$ into itself are birational
(generic existence and uniqueness of left and right divisions).}

	The regularization theorem for the algebraic pre-groups can be
stated as follows (see [30]; [34], Theorem~3.1).

\claim {Lemma~5.11}
	There exists a birational homomorphism $\tau$ between $Q$
and an algebraic group $\tilde G$.
\endclaim

	{\bf Remark.} By a {\bf birational homomorphism} we mean a birational
correspondence $\tau$ such that $\tau(uv)=\tau(u)\tau(v)$
whenever all expressions are defined (cf. [34], Definition~3.2).

\bigskip
{\it 3. Properties of the action on $V$.}
\medskip

\claim {Lemma~5.12}
	There exists a rational action $\tilde\rho\colon Q\times V \to V $,
i.e. $\tilde\rho(vw,x) \hfill\break
 = \tilde\rho(v,\tilde\rho(w,x))$ for generic choice
of $(v,w,x)\in Q\times Q\times V $
such that the following diagram is commutative whenever the mappings are
defined:
$$
\matrix
G\times V               &  \mapright{\rho}       &   V            \\
\mapdown{\phi\times id} &                        &  \mapdown{id}  \\
Q\times V               &  \mapright{\tilde\rho} &   V            \endmatrix
$$
\endclaim

	{\bf Proof.}
By Lemma~5.4, there exists a Zariski open dense subset $Q'\subset Q$
such that, for every $v\in Q'$, the fibre $F_v$ is
the closed graph of a birational mapping $\rho_v\colon V \to V $.
These birational mappings together define the action
$\tilde\rho\colon Q'\times V \to V $ which extends to a rational
mapping $\tilde\rho\colon Q\times V \to V $. The commutativity
of the diagram follows from the coincidence of the closed graph
of $\rho_g,g\in U$, with the fibre $F_{\phi(g)}$.

	The property $\tilde\rho(vw,x)=\tilde\rho(v,\tilde\rho(w,x))$
is true for $v,w,vw\in \phi(U)$ and $x\in D$. Since the set of
such $(v,w,x)$ is Zariski dense in $Q\times Q\times V $,
this is true for generic choices of $(v,w,x)$.\hfill {\bf QED}

	To  simplify the notation, we write
$\rho\colon Q\times V \to V $ instead of
$\tilde\rho\colon Q\times V \to V $. Then the property
$\tilde\rho(vw,x)=\tilde\rho(v,\tilde\rho(w,x))$ can
be written as the associativity condition $(vw)x=v(wx)$.

The ``action'' $Q\times V\to V$ is rational. Furthermore,
it may happen that an element $v\in Q$ does not
define a birational automorphism of $V$. This is not the case, however,
if we replace $Q$ by $\tilde G$.

\claim {Lemma~5.13}
	Let $\tilde G\times V \to V$ be the rational action which is
induced by the action $\rho\colon Q\times V \times V$
via the birational homomorphism $\tau\colon Q\to \tilde G$.
Then, for every element $v\in\tilde G$,
the restriction $\rho_v\colon V\to V$ is a
birational automorphism of $V$
and the product $vw$ in $\tilde G$ corresponds to the
composition $\rho_v\circ\rho_w$ of these automorphisms.
Moreover, the action of $\tilde G$ is effective.
\endclaim

	{\bf Proof.} Let $Q'\subset Q$ be the open dense subset given
by Lemma~5.4 $Q''\subset Q'$ an open dense subset where
the birational homomorphism $\tau\colon Q\to \tilde G$ is biregular.
We can regard $Q''$ as a Zariski open dense subset of $\tilde G$.

	Let $v\in \tilde G$ be arbitrary and $w\in vQ''\cap Q''$.
Then $v=wu^{-1}$ for $w,u\in Q''$. By Lemma~5.4, the fibres
$F_w$ and $F_u$ coincide with closed graphs of $\rho_w$ and $\rho_u$.
Therefore there exist points in $w\times V $ and $u\times V $
where $\rho$ is defined. By Lemma~5.8,
$\rho_v=\rho_w\circ\rho_u^{-1}$, which is also a birational automorphism
of $ V $.

	By Lemma~5.8, one has $\rho(vw,x)=\rho(v,\rho(w,x))$ for
all $(v,w,x)$ in a Zariski dense subset of
$\tilde G\times \tilde G\times  V $. Therefore this is
true for all values of $(v,w,x)$ whenever the expressions are defined.
This implies $\rho_{vw}=\rho_v\circ\rho_w$ for all $v,w\in \tilde G$.

	Assume that $\rho\colon \tilde G\times V \to V $ has a kernel $K$ and
take $k\ne 1\in K$. Let $Q''\subset Q$ be a Zariski open dense subset
where the birational homomorphism $\tau\colon Q\to \tilde G$ is biregular.
We can regard
$Q''$ as a Zariski open dense subset of $\tilde G$.
Let $v=kw\in Q''\cap kQ''$ be an arbitrary point. Since $k$ is in the kernel,
$\rho_{kw}=\rho_w$. On the other hand, $w$ and $kw$ are different points
in the Chow scheme $C\supset Q''$ with different fibres. Since the fibres
are the closed graphs of corresponding automorphisms, this is a
contradiction. \hfill {\bf QED}

	Recall that $V$ is an algebraic pre-transformation $\tilde G$-space
([34], Definition~4.1) if\par
\item {1)}{ for generic $(v,w,x)\in \tilde G\times \tilde G\times V$,
	both expressions $(vw)x$ and $v(wx)$ are defined and equal
	(generic associativity condition);}
\item {2)}{ the mapping $(v,x)\mapsto (v,vx)$ from $Q\times V$ into itself
	is birational.}

\claim {Corollary~5.1}
	$V$ is an algebraic pre-transformation $\tilde G$-space.
\endclaim

\bigskip
{\it 4.The homomorphism from $G$ into $\tilde G$.}
\medskip

	Up to now we constructed an open dense subset $U\subset G$ and
a local homomorphism $\phi\colon U\to Q$. We wish to extends $\phi$
to a homomorphism from $G$ into $\tilde G$, which is compatible with
the action on $V$.

\claim {Lemma~5.14}
	Let $Q''\subset Q$ be a Zariski open dense subset.\par
\item {1)}{
For every $g\in G$ there exist two points $v,w\in Q''$ such that
$\rho_g=\rho_v\circ\rho_w^{-1}$;}
\item{2)}{
The points $v,w$ can be chosen to
be in $\phi(U)$;}
\item {3)}{
If $g_m\to g_0$
is any convergent sequence in $G$,
the corresponding sequence $v_m,w_m\in Q''$ can be
chosen to converge to some $v_0,w_0$ with
$\rho_{g_0}=\rho_{v_0}\circ\rho_{w_0}^{-1}$.}
\endclaim

	{\bf Proof.} Let $g\in G$ be fixed and $Q'\subset Q$ be given by
Lemma~5.4. We can assume $Q''\subset Q'$.
Let $F\subset Q''\times V \times V $ be the universal family over $Q''$,
which, by Lemma~5.4, consists of closed graphs of
birational automorphisms $\rho_v\colon V \to V $.
Analogous to the formula (5.2) we can consider the family $F'$
of compositions $\rho_g\circ\rho_v,v\in Q''$.
Let $Q_0$ denote an irreducible component of $Q$.
As in the proof of Lemma~5.5, we conclude that
the fibre $F'_v$ coincides with the closed graph of $\rho_g\circ\rho_v$
for all $v$ from a Zariski open subset $Q''_g\subset Q''\cap Q_0$.
By Proposition~5.3, this family yields a rational mapping
$r_g\colon Q''_g\to C$.

	We wish to prove that $r_g(Q''_g)\subset Q$. For this we return to our
group $G$. Let $U\subset G$ be the chosen open dense subset.
Then the translation $gU$ is also an open dense subset of $G$ and
so is the intersection $U':=gU\cap U$. This implies that
$\phi(U')$ is Zariski dense in $Q$ and therefore $\phi(U')\cap Q''$
is Zariski dense in $Q''$. Now, for every
$v=\phi(h)\in \phi(U')\cap Q'',h\in U'$
the fibre $F'_v$ is the closed graph of $\rho_{gh}$. Since $gh\in U$,
one has $r_g(v)=\phi(gh)\in Q$. By the density of $\phi(U')\cap Q''_g$,
the image of $r_g(Q''_g)$ lies in $Q$.

	Since the compositions of $\rho_g$ with different automorphisms of $V$
are different, the mapping $r_g\colon Q''_g\to Q$ is injective. Therefore
the image $r_g(Q''_g)$ intersects the open dense subset $Q''$. Let
$v\in Q''\cap r_g(Q''_g)$ be an arbitrary point.
The fibre $F_v$ is the closed graph of the birational
automorphism $\rho_v$ and, at the same time, is the closed graph of
$\rho_g\circ\rho_w$, where $w\in Q''$.
This means $\rho_g=\rho_v\circ\rho_w^{-1}$ which finishes the proof
of the part~1.

	The point $v\in Q''\cap r_g(Q''_g)$ can be chosen to lie in
$\phi(U\cap gU)$. Then $v,w\in \phi(U)$ and the part~2 is also proven.

	If we are given a convergent sequence $g_m\to g_0$, we can choose
a point $v\in Q''$ such that $v\in r_{g_m}(Q''_{g_m})$
for all $m=0,1,\ldots$. Then all $w_m\in Q''_{g_m}$ lie in the component
$Q_0$ which is included in a single component of the Chow scheme.
This means that the degree of $\rho_{w_m}$ is constant.
The convergence $\rho_{g_m}\to \rho_{g_0}$ implies $\rho_{w_m} \to\rho_{w_0}$.
By Lemma~5.3, the family of closed graphs of $\rho_{w_m}$
is closed. By Proposition~5.2, $w_m\to w_0,m\to \infty$.
This proves the part~3.
\hfill {\bf QED}

\claim {Lemma~5.15}
	There exists a continuous homomorphism $\phi\colon G\to \tilde G$,
such that $\rho_{\phi(g)}=\rho_g$ for all $g\in G$.
The image $\rho(G)$ is Zariski dense in $\tilde G$.
\endclaim

	{\bf Proof.} Let $Q''\subset Q$ be a Zariski open dense subset where
the birational mapping homomorphism $\tilde G$ is biregular. We can identify
$Q''$ with a subset of $\tilde G$. By Lemma~5.14, applied to
the set $Q''\subset Q$ and an element $g\in G$, one has
$\rho_g=\rho_v\circ\rho_w^{-1}$ for $v,w\in \tilde G$.
By Lemma~5.13, $\rho_g=\rho_{vw^{-1}}$.
Then we define $\phi(g):= vw^{-1}$. Since, by Lemma~5.13,
the action of $\tilde G$ is effective, this definition of $\phi(g)$ is
independent of the choices of $v$ and $w$.

	The property $\rho_{\phi(g)}=\rho_g$ is satisfied by construction
of $\phi$. By Lemma~5.13, $\phi$ is a homomorphism.
The continuity of $\phi$ follows from Lemma~5.14, part~3.

	The image $\phi(G)$ contains the image $\phi(U)$, which is Zariski
dense in $Q$. If $Q''\subset Q$ is a Zariski open dense subset, where
the isomorphism between $Q$ and $\tilde G$ is biregular, the intersection
$Q''\cap \phi(U)$ is Zariski dense in $Q$. The set $Q''$ can be regarded
as a Zariski open dense subset of $\tilde G$. This which yields the density
of $Q''\cap \phi(U)$ and therefore of $\phi(G)$. \hfill {\bf QED}

\bigskip
{\it 5. The regularization of the action $\tilde G\times V\to V$.}
\medskip

	Let  $D\subset V$ be as in Theorem~$3'$.
We noted in Corollary~5.1 that $V$ is an algebraic pre-transformation
$\tilde G$-space. The theory of A.~Weil (see Theorem~4.1 in [34]) gives the
existence of
the regularizations of algebraic pre-transformation spaces which
are regular at the so-called points of regularity.
Recall that a point $x$ in an algebraic pre-transformation $\tilde G$-space
$V$ is called a {\bf point of regularity} if the mapping $x'\mapsto ux'$
from $V$ into itself is biregular at $x'=x$ for generic $u\in \tilde G$
(see [34], Definition~4.3).

	If $v\in\phi(G),x\in D$, then the mapping $x'\mapsto vx'$ is biregular
at $x$. Since, by Lemma~5.15, $\phi(G)$ is Zariski dense in $\tilde G$,
$D$ consists of points of regularity. By Theorem~4.1 in [34],
there exists a birational regularization $\psi\colon V\to X$,
i.e. $\tilde G$ acts regularly on $X$, the mapping $\psi$
is birational on $V$, biregular on $D$ and $\tilde G$-equivariant.
In particular, $\psi|_D$ is $G$-equivariant.
This is exactly the conclusion of Theorem~$3'$.
\hfill {\bf QED}

\newsection 6. Algebraic extensions for the  case of finitely many connected
components \par

	In this section we prove Theorem~1.
Let $G$ be a Lie group of birationally extendible
automorphisms of $D\subset V$ with finitely many connected components.
For fixed $g\in G$ we obtain an $n$-dimensional
subvariety $Z_g=\nu(\Gamma_g)\subset\P_k$ which corresponds to a point
$\rho(g)$ in the cycle space $\c(\P_k)$ ([1,5]).
In order to apply the universality of the cycle space we embed
our family $Z_g,g\in G$ in a meromorphic family $Z_g, g\in \tilde G$.

	As a real analytic manifold, $G$ can be embedded totally really and
closed into a complex manifold $G'$  with $\dim_\R G = \dim_\C G'$
([4]). We wish to extend the action $G\times D\to D$
to a meromorphic mapping $\tilde G\times V\to V$, where
$\tilde G$ is a neighborhood of $G$ in $G'$.
Since $G$ is embedded totally really, the meromorpic extension is unique.
Therefore it only must be constructed locally with respect to $G$.
For the proof we utilize the following result of Kazaryan ([15]).
A subset $E\subset D'$ is called {\it nonpluripolar} if
there are no plurisubharmonic functions $f\colon D'\to \R\cup\{-\infty\}$
such that $f|_E\equiv -\infty$.

\claim {Proposition~6.1}
	 Let  $D'$  be  a domain in $\C^n$  and let $E\subset D'$ be a
nonpluripolar
subset. Let $D''$ be an open set in a
complex  manifold $X$. If $f$ is a meromorphic function on $D'\times D''$
such that $f(g,\cdot)$ extends to a meromorphic function on $X$
for  all $g\in E$, then $f$ extends to a meromorphic function in
a neighborhood of $E\times X\subset D'\times X$.
\endclaim

	We wish to prove the required extension at a point $g_0\in G$.
For this we fix a coordinate neighborhood $E\in G'$ of $g_0$
regarded as a neighborhood in $\C^p$, such that $G\cap E=\R^p\cap E=:E_{\R}$.
The map $\mu=\nu\circ(id\times \phi) \colon G\times V\to \P^k$
is real analytic on $E_{\R}\times D$ and extends therefore
to a holomorphic map in a neighborhood  $D'\times D''$
of  $E_{\R}\times D\subset \C^p\times V$.
(Here we must replace $D$ by a bit smaller neighborhood $D''\subset D$).

	The set $E_{\R}$, being an open subset of $\R^p$ , is nonpluripolar.
We apply Proposition~6.1 to the coordinates of the map $\mu$
in any affine coordinate chart in $\P_k$.
We conclude that $\mu$ extends to a meromorphic map
$\tilde\mu$ defined in a neighborhood of
$\{x_0\}\times V\subset G'\times V$ into $V$.
Since $V$ is compact, we can choose this neighborhood of the form
$\tilde G\times V$.

	Now we can apply Proposition~5.1 to the meromorphic family
$Z_g, g\in \tilde G$. We obtain a meromorphic mapping
$\phi\colon\tilde G\to {C}_n(\P_k)$. Since the number of components of $G$
is finite, $\tilde G$ can be also assumed to possess this property.
Then the image $\phi(\tilde G)$ lies also in finitely many components
of ${C}_n(\P_k)$ which means the boundness of the degree for all
$Z_g$ with $g$ in an open dense subset $U\subset \tilde G$.
By Lemma~5.1, the degree is globally bounded.
Now the application of Theorem~3 yields the algebraic extension required by
Theorem~1.  \hfill {\bf QED}

\newsection 7. The proof of Theorem~4 \par

{$1\Longrightarrow 2$.} The proof is trivial. \hfill {\bf QED}
\bigskip

{$2\Longrightarrow 3$.}
	By Theorem~3, it is sufficient to prove the boundness of
the degree. Let $G$ be a subgroup of
a Nash group $\tilde G$ such that the action
$G\times D\to D$ extends to a Nash action $\tilde G\times D\to D$.
We prove the statement for arbitrary Nash manifold $\tilde G$ and
Nash map $\tilde G\times D\to D$
(which is holomorphic for every fixed $g\in G$)
by induction on $\dim \tilde G$. This is obvious for $\dim \tilde G=0$.

	Let $U\subset \tilde G$ and $W\subset \P_k$ be Nash coordinate charts and
$\phi_j(g)\colon D\to \C$ be the $j$th coordinate in $W$ of
$\nu\circ({id}\times\phi_g)\colon D\to \P_k$
for $g\in U$ (taken on its set of definition).
Since the map $\phi_j\colon U\times D\to \C$ is Nash, it
satisfies a polynomial equation
$P_j(g,x,\phi_j(g,x))\equiv 0$ of degree $d$.
This yields nontrivial polynomial
equations of degree not larger than $d$ for all $\phi_j(g)\colon D\to \C$,
$g\in U$, outside a proper algebraic subvariety $N$. The calculation of the
required degree, i.e. the intersection number with a linear projective
subspace $L$ of codimension $n$, yields additional linear equations for
the coordinates in $W$. For $L$ generic and $g$ in the complement of another
proper subvariety $N'$, this intersection number in finite. Since the degrees
of polynomial equations for this intersection are bounded, the intersection
number is also bounded (Bezout theorem). This proves the statement
for $\tilde G=U\backslash (N\cup N')$.

	The intersection $U\cap (N\cup N')$
admits a finite stratification in lower dimensional Nash manifolds
(see e.g. [2]). By induction, the required degree is bounded
for every stratum. This proves the boundness of degree for $g\in U$.
Since the Nash atlas is finite, we obtain the required boundness
for the Nash manifold $\tilde G$. \hfill {\bf QED}

	\bigskip
	{$3\Longrightarrow 1$.} Let $\rho\colon\tilde G\times X\to X$
be an algebraic extension. We identify the open subset $D$ with its
embedding in $X$. Then we define $\tilde G'$ to be the subgroup
of $\tilde G$ which consists of all elements which leave $D$ invariant.
In general, this is not an algebraic subgroup. For our statement, it
is sufficient to prove that $\tilde G'$ is a Nash subgroup.

	We utilize the following property of semialgebraic sets
([33], Lemma~6.2).

\claim {Lemma~7.1}
	Let $A$, $B$ and $C,C'\subset A\times B$ be semialgebraic sets.
Then the set of $a\in A$ such that $C_a\subset C'_a$ is semialgebraic.
\endclaim

	Here $C_a$ and $C'_a$ denote the fibres $\{b\in B\mid (a,b)\in C\}$
and $\{b\in B\mid (a,b)\in C'\}$ respectively.
Now we set $A:=\tilde G$, $B:=D$,
$C:=(pr_{\tilde G}\times\rho)(\tilde G\times D)$, $C':=\tilde G\times D$
in Lemma~7.1. The set $C$ is semialgebraic by the Tarski-Seidenberg
Theorem ([2], Theorem~2.7.1). By Lemma~7.1, the set
$G_1:=\{g\in \tilde G\mid g(D)\subset D\}$ is semialgebraic.
Again, Lemma~7.1, applied to $A:=G_1$, $B:=D$, $C:=G_1\times D$ and
$C':=(pr_{G_1}\times\rho)(G_1\times D)$, shows that the subgroup
$\tilde G'$ is semialgebraic. Therefore it is a Nash subgroup and
the statement is proven. \hfill {\bf QED}


\references

\refer[1] D.~Barlet,
{\it Espace cycles analytique complexes de dimension finie,}
Seminare F. Norguet, Lecture Notes in Math., Springer,
{\bf 482} (1975), 1--158. \par

\refer[2] R. Benedetti and J.-J. Risler,
{\it Real algebraic and semi-algebraic sets,}
Actualites Mathematiques. Hermann Editeurs des Sciences et des Arts,
(1990).\par

\refer[3] E. Bishop,
{\it Conditions for the analyticity of certain sets,}
Mich. Math. J., {\bf 11} (1964,) 289--304.\par

\refer[4] F. Bruhat and H. Whitney,
{\it Quelques propri\'et\'e fondamentales des ensembles analytiques r\'eels,}
Comment. Helv., {\bf 33} (1959), 132--160.\par

\refer[5] F. Campana and Th. Peternell,
{\it Cycle spaces,}
In Encyclopedia of Mathematical Sciences,
{\it Several Complex Variables VII,} Springer,
{\bf 74} (1994), 319--349.\par

\refer[6] H. Cartan,
{\it Sur les groupes de transformations analytiques,}
Act.~Sc.~et~Int., Hermann, Paris, (1935). \par

\refer[7] K. Diederich,
Private communications, (1995).\par

\refer[8] R. E. Green and S. G. Krantz,
{\it Charakterization of certain weakly pseudoconvex domains with noncompact
automorphism groups,}
Lecture Notes in Math.,  Springer, {\bf 1268} (1987), 121-157.\par

\refer[9] Ph. Griffiths and J. Harris,
{\it Principles of Algebraic Geometry,}
John Wiley \& Sons, (1978).\par

\refer[10] J. Harris,
{\it Algebraic Geometry: a first course,}
Graduate Text in Math., Springer, {\bf 133} (1993).\par

\refer[11] R. Hartshorne,
{\it Algebraic ge\-ometry,}
Graduate Texts in Mathe\-matics, \hfill\break
Sprin\-ger, {\bf 52} (1977).\par

\refer[12] P. Heinzner,
{\it Geometric invariant theory on stein spaces,}
Math. Ann., {\bf 289} (1991), 631--662.\par

\refer[13] P. Heinzner and A. Ianuzzi,
{\it Integration of local actions on holomorphic fiber spaces,}
preprint, (1995).\par

\refer[14] W. Kaup, Y. Matsushima, and T. Ochiai,
{\it On the automorphisms and equivalenes of generalized siegel domains,}
Amer. J. Math., {\bf 92} (1970), 475--498.\par

\refer[15] M. V. Kazaryan,
{\it Meromorphic continuation with respect to groups of variables,}
Math. USSR-Sb., {\bf 53} (1986), 385--398.\par

\refer[16] F. Knop, H. Kraft, D. Luna, and T. Vust,
{\it Local properties of algebraic group actions,}
In H. Kraft, P. Slodowy, and Tonny A. Springer (editors),
{\it Algebraic Transformation Groups and Invariant Theory,}
DMV-Seminar, Birkh\"auser,
{\bf 13} (1989), 63--76.\par

\refer[17] A. Kodama,
{\it On the structure of a bounded domain with a special boundary point (II),}
Osaka~J.~Math., {\bf 24} (1987), 499--519.\par

\refer[18] J. J. Madden and C. M. Stanton,
{\it One-dimensional {N}ash groups,}
Pacific Journal of Math., {\bf 154(2)} (1992), 341--344.\par

\refer[19] R. Narasimhan,
{\it Several complex variables,}
Chicago Lectures in Mathematics, Univ. of Chicago Press, (1971).\par

\refer[20] R. Penny,
{\it The structure of rational homogeneous domains in $\C^n$,}
Ann. Math., {\bf 126} (1987), 389--414. \par

\refer[21] S. Pinchuk,
{\it The scaling method and holomorphic mappings,}
Proc. Symp. Pure Math., in E. Bedford, J. P. d'Angelo et al (Ed.),
Several Complex Variables and Complex Geometry,
{\bf 52(1)} (1991), 151--161.\par

\refer[22] J. P. Rosay,
{\it Sur une caract\'erisation de la boule parmi les domaines de $\C^n$
par son groupe d'automorphismes,}
Ann. Inst. Fourier, {\bf 29} (1979), 91--97. \par

\refer[23] M. Rosenlicht,
{\it Some basic theorems on algebraic groups,}
Amer. J. Math., {\bf 78} (1956), 401--443,\par

\refer[24] O. Rothaus,
{\it The construction of homogeneous convex cones,}
Bull. Amer. Math. Soc., {\bf 69} (1963), 248--250,\par

\refer[25] I. R. Shafarevich,
{\it Basic algebraic geometry,} Springer, (1974).    \par

\refer[26] M. Shiota,
{\it Nash manifolds,}
Lect. Notes Math., Springer, {\bf 1269} (1987).  \par

\refer[27] H. Sumihiro,
{\it Equivariant completion {II},}
Math. Kyoto Univ., {\bf 15} (1975), 573--605.\par

\refer[28] E. B. Vinberg, S. G. Gindikin, and I. I. Pyatetskii-Shapiro,
{\it Classification and canonical realization of complex bounded homogeneous
domains,}
Trans. Moscow. Math. Soc., {\bf 12} (1963), 404--437.\par

\refer[29] S. Webster,
{\it On the mapping problem for algebraic real hypersurfaces,}
Inventiones math., {\bf 43} (1977), 53--68.\par

\refer[30] A. Weil,
{\it On algebraic group of transformations,}
Amer. J. of Math., {\bf 77} (1955), 355--391.\par

\refer[31] B. Wong,
{\it Characterization of the unit ball in $\C^n$ by its automorphism group,}
Invent. Math., {\bf 41} (1977), 253--257.\par

\refer[32] B. Wong,
{\it Charakterisation of the bidisc by its automorphism group,}
Amer. J. of Math., {\bf 117(2)} (1995) 279--288. \par

\refer[33] D. Zaitsev,
{\it On the automorphism groups of algebraic bounded domains,}
Math. Ann. {\bf 302} (1995), 105-129.\par

\refer[34] D. Zaitsev,
{\it Regularizations of birational group operations in sense of Weil,}
Preprint, (1995).\par


\endreferences

\enddocument